\documentclass[11pt]{amsart}

\usepackage{fullpage,url,amssymb,enumerate,colonequals}
\usepackage{graphicx} % Required for inserting images
\usepackage{amsmath}
\usepackage{amsthm}
\usepackage{amssymb}
\usepackage{tikz-cd}
\usepackage{tikz}
\usepackage{tabularx}
\usepackage{footnote}
\usepackage{hyperref}
\usepackage{lineno}
\usepackage{xcolor}
\usepackage[dvipsnames]{xcolor}
\makesavenoteenv{table}
\makesavenoteenv{tabularx}

\DeclareMathOperator{\Gal}{Gal}

\DeclareMathOperator{\divv}{div}
\DeclareMathOperator{\Aut}{Aut}
\newtheorem{theorem}{Theorem}
\newtheorem{proposition}[theorem]{Proposition}

\theoremstyle{theorem}

\theoremstyle{definition}
\newtheorem{remark}{Remark}

\newtheorem*{Ack}{Acknowledgements}

\newcommand{\Q}{\mathbb{Q}}

\title{Three-Torsion Subgroups and Wild Conductor Exponents of Plane Quartics}
\author{Elvira Lupoian and James Rawson}
\date{}
\address{Elvira Lupoian \\ Department of Mathematics, University College London, London \\ United Kingdom}
\email{e.lupoian@ucl.ac.uk}
\address{James Rawson \\ Mathematics Institute, University of Warwick, Coventry \\ United Kingdom}
\email{james.rawson@warwick.ac.uk}

\begin{document}
\begin{abstract}
In this paper we give an algorithm to find the 3-torsion subgroup of the Jacobian of a smooth plane quartic curve with a marked rational point. We describe $3-$torsion points  in terms of cubics which triply intersect the curve, and use this to define a system of equations whose solution set corresponds to the coefficients of these cubics. 
We compute the points of this zero-dimensional, degree $728$ scheme first by approximation, using homotopy continuation and Newton-Raphson, and then using continued fractions to obtain accurate expressions for these points. 
We describe how the Galois structure of the field of definition of the $3$-torsion subgroup can be used to compute local wild conductor exponents, including at $p=2$. 
\end{abstract}
\maketitle

\section{Introduction}
The celebrated result of Mazur \cite{mazur1977modular} gives a complete classification of rational torsion subgroups of elliptic curves. The general study of torsion of abelian varieties is considered out of reach. A number of works have aimed to construct Jacobians of curves with rational torsion points of large order, see for instance \cite{bernard2009jacobians}, \cite{howe2000large}, \cite{flynn1990large}, \cite{leprevost1992torsion}, \cite{howe2015genus}, and there has been some work towards computing the entire rational torsion subgroup of Jacobians of hyperelliptic curves, see \cite{stoll1998height} and \cite{j2023computing}. 
An interesting variant of this question is the study of $n$-torsion points, for a fixed natural number $n$. 
For Jacobians of genus 1 curves, this is a largely elementary problem, since explicit equations are known for the Jacobian (as it is an elliptic curve) and the group law is accessible. For Jacobians of higher genus curves, the problem becomes much more difficult. Perhaps the clearest description is that of $2$-torsion, since for the hyperelliptic curves, it can be easily deduced from the Weierstrass model, and for non-hyperelliptic curves $2$-torsion points correspond to multi-tangent hyperplanes to the curve, as described in \cite{dog}.  Torsion points of higher order on genus $2$ are explicitly described by Flynn, Testa and Bruin \cite{bruin2014descent} and Flynn \cite{flynn2015descent}.
 The description of $3$-torsion is generalised in \cite{lupoianthree}, for hyperelliptic curves of genus $3$. For hyperelliptic curves, such description of torsion can be obtained using a Weierstrass model for the curve and careful analysis of Riemann-Roch spaces. However, when the curve is not hyperelliptic, the description of torsion points tends to become more directly geometric. Explicitly computing such torsion subgroups has important applications. For instance, it is used for descent purposes in \cite{bruin2014descent} and \cite{BS2tors}; and also for conductor computations in \cite{dokdor}.

In this paper we study the $3$-torsion subgroup of genus $3$, non-hyperelliptic curves. We begin by showing that for plane quartics over a number field $K$, which have a $K$-rational point, after a suitable normalisation, $3$-torsion points on the Jacobian correspond to cubics of the form 
\begin{center}
    $ \alpha_{1}x^3 + \alpha_{2}z^{3} + \alpha_{C}\alpha_{6} x^2y + \alpha_{3}x^2z + \alpha_{4}xz^2 + \alpha_{5}xyz + \alpha_{6}y^2z + \alpha_{7}yz^2 
    $
\end{center}
which intersect the curve with multiplicity divisible by $3$ at each point and $\alpha_{C}$ is a constant depending explicitly on the canonical model of the curve. Computing the $3$-torsion subgroup of the Jacobian reduces to computing the coefficients of such cubics. To achieve this, we derive equations whose solutions parametrise the coefficients of such cubics, which we solve using a two step approach. First we compute high precision complex approximations of such solutions, and we obtain precise algebraic expressions using lattice reduction or continued fractions.  This two step approach was previously used for computing torsion points in \cite{lupoiantwo} and \cite{lupoianthree}.

As in \cite{dokdor} and \cite{lupoianthree}, we make use of this explicit computation of $3$-torsion to determine the wild conductor exponents at $2$. Recall that for any smooth and projective curve over $\Q$ the conductor is a representation theoretic quantity attached to a curve, defined as a product over the primes of bad reduction of the curve. The problem of computing the conductor, reduces to computing the exponent, the \textit{local conductor exponent}, of each prime of bad reduction. 

In the case of elliptic curves, Tate's algorithm \cite{silv2} 
can be used to compute the conductor explicitly. For hyperelliptic curves, most of the conductor can be determined using the theory of cluster pictures \cite{ddmm}. Moreover, in the case of genus $2$ curves, the conductor can be completely determined using the algorithm of Dokchitser-Doris \cite{dokdor}. 

Computationally, the conductor exponent at any prime $p$ can be viewed as the sum of two parts
\begin{center}
    $n_{p} = n_{\text{tame}, p} + n_{\text{wild}, p}$
\end{center}
the tame part $n_{\text{tame}, p}$, which can be determined explicitly from a regular model of the curve, and the wild part $n_{\text{wild}, p}$, which can be determined explicitly from knowing the Galois action on $J[l]$ for any prime $l \ne p$.  Thus in the case of plane quartics, our explicit $3$-torsion calculations can be used to determine $n_{\text{wild}, 2}$. Furthermore, as described in the final section of this paper, the classical description of $2$-torsion points in terms of the bitangents to the curve, can be used to compute $n_{\text{wild}, p}$ for all $p \ne 2$. 

The code referenced throughout this paper can be found in the following repository:
\begin{center}
\href{https://github.com/ElviraLupoian/Genus3Conductors}{https://github.com/ElviraLupoian/Genus3Conductors}
\end{center}

\subsection*{Outline}
This paper is organised as follows. In Section \ref{3tors} we give a description of $3$-torsion points using sections of the canonical class and a marked rational point. In Section \ref{3torscomps} we describe a method for computing $3$-torsion points. We derive a system of equations whose solutions correspond to $3$-torsion points, which we solve by first approximating the solutions, as tuples of complex numbers, and then by using such approximations and lattice reduction/continued fractions to obtain algebraic expressions for our torsion points. We give an overview of some examples in Section \ref{egs}.
In Section \ref{galrep} we describe how the local Galois representation of $J[3]$ can be computed from our data. 
In Section\ref{2torsion} we give a brief overview of the description of the $2$-torsion subgroup, and explain how one compute the local Galois representation at $2$.
In the final section we explain how our calculations can be used to determine the wild conductor exponents at all primes. 

\subsection*{Conventions and Notation} 
For the rest of this paper, unless stated otherwise, we write 
$C$ for a smooth and projective, non-hyperelliptic curve of genus $3$, defined over a number field $K$. The curve is embedded in $\mathbb{P}^{2}$ by the canonical embedding and it's described by a single degree $4$, homogeneous polynomial $f \in K[x,y,z]$. We write $J$ for the Jacobian of $C$. Recall that this is a $3$-dimensional abelian variety over $K$, whose points correspond to linear equivalence classes of degree $0$ divisors on the curve, by the Picard functionality.

\begin{Ack}
 The authors are grateful to Raymond van Bommel, Tim Dokchitser, Samir Siksek and Damiano Testa  for helpful conversations. We thank the anonymous referees for their careful reading of this paper. Their comments and suggestions have improved this work greatly. Both authors are grateful for the financial support provided by the UK Engineering and Physical Sciences Research Council. The first named author is supported by the EPSRC Doctoral Prize fellowship EP/W524335/1. The second named  author is supported by the  EPSRC studentship EP/W523793/1.
 \end{Ack}

\section{Three-Torsion Points and Sections of $K_{C}$} \label{3tors}
The group of $2$-torsion points is classically generated by differences of odd theta characteristics of the curve, which in turn correspond to bitangent lines to the curve, as described in Section~\ref{2torsion}. The problem of describing $3$-torsion points on $J$ is not as straightforward. Notably, a `nice geometric description' analogous to bitangents for $3$-torsion points does not exist, since $\text{GSp}_{6}(\mathbb{F}_{3})$ has no small index subgroups, with the smallest having order $364$, and corresponding to `projective` $3$-torsion points. In this section we give an elementary description of $3$-torsion points, using sections of the canonical divisor. 

From now on, we additionally assume that $C$ has a $K$-rational point $P_{0}$. By a suitable change of coordinates we assume $P_{0} = ( 0 : 1 :0)$ and that the tangent line to $C$ at $P_{0}$ is given by $z=0$. Recall that we work with the canonical embedding of $C$ into $\mathbb{P}^{2}$, and suppose that the curve is cut out by a quartic $f \in K[x,y,z]$, 
\begin{center}
    $C: f(x,y,z) = a_{0}x^4 + a_{1}x^3y + a_{2}x^3z + a_{3}x^2y^2 + a_{4}x^2yz  + a_{5}x^2z^2 + a_{6}xy^3 + a_{7}xy^2z + a_{8}xyz^2 + a_{9}xz^3 + a_{10}y^4 + a_{11}y^3z + a_{12}y^2 z^2 +a_{13}yz^3 + a_{14}z^4= 0.$
\end{center}
As $P_{0}= (0:1:0) \in C(K)$ and $z=0$ is the tangent line at this point, it is necessary that $a_{10} = a_6 = 0$ and $a_{11} \ne 0$ else $C$ would be hyperelliptic. For the rest of this section we write $\alpha_{f} := \frac{a_{3}}{a_{11}}.$

\begin{theorem} \label{desc}
Let $C$ be as above. Three-torsion points of the Jacobian of $C$ correspond to cubics of the form 
  \begin{center}
    $ \alpha_{1}x^3 + \alpha_{2}z^{3} + \alpha_{f} \alpha_{6} x^2y + \alpha_{3}x^2z + \alpha_{4}xz^2 + \alpha_{5}xyz + \alpha_{6}y^2z + \alpha_{7}yz^2 
    $
  \end{center}
  which intersect $C$ with multiplicity divisible by $3$ at each point. 
\end{theorem}

\begin{proof}
 Fix an effective representative $K_C$ of the canonical divisor whose support contains $P_{0}$. For any $D \in J[3] \setminus \{ 0 
 \}$, we note that $D + K_C$ is a divisor of degree $4$ and thus by the Riemann-Roch theorem, there exists an effective, degree $4$ divisor $ D_{0}$, whose support contains $P_{0}$ and such that $D + K_C \sim D_{0}$. Thus, there exists a non-constant function on the curve $h$, satisfying
 \begin{center}
     $\divv h = 3D_{0} - 3K_C$.
 \end{center}
 As $D_{0}$ and $K_C$ are assumed to be effective, $h \in \mathcal{L} (3K_C)$ and hence, it is precisely a cubic. As this cubic passes through $P_{0}$, the coefficient of $y^{3}$ is necessarily zero. Moreover, as $z=0$ is the tangent line at this point, we have the additional condition that the coefficient of $y^2x$ is also zero. Lastly, the fact that the  cubic intersect the quartic with multiplicity $3$ gives us the relation between the coefficient of $x^2y$ and $y^2z$. Working in the affine patch obtained by dehomogenising with respect to $y$ shows that at $P_0$, $x$ is a uniformiser and $z =  \alpha_f x^2 \mod x^3$, and for the cubic to vanish to order 3 at $P_0$, this equality must hold there too.
 \end{proof}

We must note that not all cubics of this form give rise to non-trivial torsion elements. With notation as in the above proof, if $D_{0} - K_C \sim 0$, then $D_{0}$ is a canonical divisor and so the intersection of $C$ with a line. Since it contains $P_{0}$, the only lines through this point are the lines of constant $x/z$, hence cubics of the form $(x - az)^3$ are excluded from the above correspondence.

Moreover, representations of $3$-torsion points by cubics of the above form are almost always unique. If $D_{0} - P_{0}$ is not a special divisor (that is, $l(K_C - (D_0 - P_0)) = 0$) and $D_0 - K_C \sim D_1 - K_C$, then $D_0 - P_{0} \sim D_1 - P_{0}$, and so there is equality of divisors as neither moves. As the divisors are the same, the cubics cutting them out must be the same up to scaling. If $D_{0} - P_{0}$ is a special divisor, then it can be written as $K' - Q$, for some canonical divisor $K'$ and point $Q$. It follows that, up to linear equivalence, $D \sim P_{0} - Q$, and so $3D \sim 3P_{0} - 3Q \sim 0$. This shows that both $P_0$ and $Q$ are flexes. The tangent at $Q$ meets $C$ at another point, $R$, and the divisors equivalent to $3Q$ are the intersection of $C$ with lines through $R$. As $3P_{0} \sim 3Q$, the tangent to $P_{0}$ must meet $C$ again at $R$. So as long as $C$ has either a rational point which is not a flex, or a rational flex whose flex line does not meet another in this way, this method will compute all non-trivial three torsion, and find unique expressions for each element.

\begin{remark}
We note that this description of torsion can be extended to higher torsion orders and genera: fix a rational degree $g - 2$ effective divisor, $\delta$, an effective representative of the canonical, $K$, such that $\delta$ is in the support of $K$, then for any $n$-torsion divisor, $D$ associate to it an effective divisor $D_0$ such that $D + K \sim D_0$ and $\delta$ is in the support of $D_0$. These divisors $D_0$ can then be identified with the sections of $nK$ whose vanishing orders are multiples of $n$. For plane degree $d$ curves, this is equivalent to intersecting the curve with degree $n(d - 3)$ curves, and for a canonically embedded curve this is intersections with degree $n$ hyperplanes.
\end{remark}

\section{Computing Three Torsion Points} \label{3torscomps}
In this section we give an overview of method for computing cubics as in Theorem \ref{desc} and hence the $3$-torsion subgroup of the Jacobian. Our method is as follows. We first derive a system of equations whose solutions parametrise the coefficients of such cubics. Theoretically, Gr\"obner basis techniques can be used to solve the system of equations outlined above. However, due to the degree of the scheme, this methodology is impractical. Instead, we employ a two step process. First, we compute high precision complex approximations, using homotopy continuation and Newton-Raphson. Then, we obtain algebraic expressions for the points via continued fraction methods or by lattice reduction.
\subsection{Schemes of $3$-torsion points}
Computing the $3$-torsion subgroup of $J$ reduces to computing coefficients of the cubics described in Theorem \ref{desc}. Our first step is to define a system of equations, whose solutions are precisely these coefficients. 

We work on the affine patch $\{ z =1 \}$ and set $F(x,y) = f (x, y,1)$. Recall that we wish to find cubics of the form 
\begin{center}
  $ \alpha_{1}x^3 + \alpha_{2} + \alpha_{6} \alpha_{f} x^2y + \alpha_{3}x^2 + \alpha_{4}x + \alpha_{5}xy + \alpha_{6}y^2 + \alpha_{7}y 
    $
\end{center}
which intersect the curve in three points, each with multiplicity $3$. Here $ \alpha_{f}$ is the previously defined constant, which is explicitly determined by the coefficients of the defining equation.  Suppose that $\alpha_{6} \ne 0$ and rescale our cubic to one of the form 
\begin{center}
  $ y^{2} - ( \beta_{1}x^3 + \beta_{2} + \beta_{f} x^2y + \beta_{3}x^2 + \beta_{4}x + \beta_{5}xy + \beta_{6}y ) = y^{2} - s_{1}(x,y)$.
\end{center}
We replace all the $y^2$ terms in the expression of $F(x,y)$ by $s_{1}(x,y)$, noting that in the resulting expression, all monomials are divisible by at most $y^2$. Repeating this procedure, the resulting equation is linear in $y$ and of the form 
\begin{center}
    $F(x,y) = s_{2}(x)y + s_{3}(x)$. 
\end{center}
Substituting $y = - s_{3}(x)/s_{2}(x)$ in the equation defining the cubic gives a rational function whose numerator is a degree $9$ polynomial $s_{4}(x)$. If the cubic intersects the curve in $3$ points, each occurring with multiplicity $3$, then $s_{4}$ is necessarily a cube, and we set 
\begin{center}
    $s_{4}(x) = \beta_{7}(x^3 + \beta_{8}x^2 + \beta_{9}x + \beta_{10})^{3}$.
\end{center}
Equating coefficients in the above equations gives 10 equations $e_{1}, \ldots, e_{10}$ in $\beta_{1}, \ldots, \beta_{10}$, and we note that $\beta_{1} \ldots \beta_{6}$ define the desired cubic. To avoid possible singularities in our system, we rule out the possibility that $s_{2}$ and $s_{3}$ have common factors by imposing the condition that their resultant is non-zero. To achieve this we add the following simple equation to our scheme 
\begin{center}
    $\beta_{11} \text{Res} (s_{2}, s_{3}) + 1 =0$.
\end{center}
where $\text{Res}(s_{2}, s_{3})$ denotes the resultant of the two polynomials.

\begin{remark}
    The polynomial $s_4$ could also be computed as the resultant (over $y$) of $f$ and $F$, however, this would not let us compute this additional resultant condition.
\end{remark}
\begin{remark}
This system of equations is not guaranteed to parametrise all $3$-torsion points, for instance this can occur if $\alpha_{6} = 0$ for a particular cubic. However, we note that for our calculations the $3$-torsion points obtained from this system are almost always sufficient, as we explain later. 
\end{remark}
 \subsection{Complex Approximations}
Approximate complex solutions of a system as above can be 
computed by homotopy continuation; a method of obtaining approximate solutions to a system of equations by deforming the known solutions of a similar system. A detailed explanation of this theory can be found in \cite{HC}. Homotopy continuation is efficiently implemented in the \texttt{Julia} package HomotopyContinuation.jl (see \cite{JulHC}), and this implementation gives approximates solutions which are accurate to 16 decimal places. We improve the precision of these approximations by the classical method of Newton-Raphson, see \cite[Page 298]{stoerNA} for a comprehensive overview of this. 

\subsection{Lattice Reduction and Precise Expressions}
When the corresponding cubic is defined over a field of small degree, we use lattice reduction to determine the minimal polynomials of the coefficients. This is a standard method, see \cite[Chapter 6]{smart} for details. We give a brief overview  of this here.

Let $x$ be an algebraic number and suppose that it is approximated to high precision by a complex number $\alpha$. We want to use $\alpha$ to find $n \in \mathbb{N}$  and $c_{0}, \ldots, c_{n} \in \mathbb{Z}$ such that 
\begin{center}
    $c_{n}x^n + \ldots + c_{1}x + c_{0} =0$.
\end{center}

Assume $\alpha$ is accurate to $k$ decimal places, and fix $k' < k$ such that $|[10^{k'} \alpha^i] - 10^{k'} x^i| \leq 1$, where $[]$ denotes the closest (Gaussian) integer. If $10^{k'}$ is much larger than the coefficients $c_i$, then the shortest vector in the lattice spanned by the columns of the following matrix should be $(c_n, \ldots, c_1, \mathrm{Re}(\gamma), \mathrm{Im}(\gamma) )$ where $\gamma = c_n [10^{k'} \alpha^n] + \ldots + c_1 [10^{k'} \alpha] + c_0 10^{k'}$.
$$\begin{pmatrix}
    1 & \ldots & 0 & 0 \\
    0 & \ldots & 0 & 0 \\
    \vdots & \ddots & \vdots & \vdots \\
    0 & \ldots & 1 & 0 \\
    [10^{k'} \mathrm{Re}(\alpha^n)] & \ldots & [10^{k'} \mathrm{Re}(\alpha)] & 10^{k'} \\
    [10^{k'} \mathrm{Im}(\alpha^n)] & \ldots & [10^{k'} \mathrm{Im}(\alpha)] & 0
\end{pmatrix}$$
This should heuristically be the shortest vector since in any other vector, the contributions from the last row will dwarf the contributions from the earlier rows. 
\begin{remark}
    When the imaginary part is zero, this last row is unnecessary and can be ignored.
\end{remark}

In practice, neither $n$ nor the sizes of the coefficients are known. We therefore work upwards from $n = 2$, and making a suitably large choice of $k'$, until we find a combination of $n$ and $k'$ such that the error term, $\gamma$, is disproportionately small.

Another use case of lattice reduction techniques is when $x$ (as before) is known to lie in a field $K$ with a given embedding into $\mathbb{C}$, and to express $x$ in terms of a basis for $K$, using $\alpha$. A particular use case for us will be as follows: if the minimal polynomial for some $\beta_1$ is known, how can $\beta_2$ be expressed in terms of the power basis for $\mathbb{Q}(\beta_1)$?

We proceed essentially as before. Let $1, x_1, \ldots, x_n$ be a basis for $K$, and $\alpha$ be a complex approximation to $x$, where $c_{n + 1} x = c_0 + c_1 x_1 + \ldots + c_n x_n$, with $c_i \in \mathbb{Z}$, and $\alpha$ is accurate to $k$ decimal places. If $k' < k$ is such that $10^{k'}$ is much larger than the $c_i$, then $(-c_{n + 1}, c_n, \ldots, c_1, \mathrm{Re}(\gamma), \mathrm{Im}(\gamma))$ will be the shortest vector in the lattice spanned by the columns of the following matrix, where $\gamma = -c_{n + 1}[10^{k'} \alpha] + c_n [10^{k'} x_n] + \ldots + c_1 [10^{k'} x_1] + c_0 [10^{k'}]$.
$$\begin{pmatrix}
    1 & 0 & \ldots & 0 & 0 \\
    0 & 1 & \ldots & 0 & 0 \\
    \vdots & \vdots & \ddots & 0 & 0 \\
    0 & 0 & \ldots & 1 & 0 \\
    [10^{k'} \mathrm{Re}(\alpha)] & [10^{k'} \mathrm{Re}(x_n)] & \ldots & [10^{k'} \mathrm{Re}(x_1)] & 10^{k'} \\
    [10^{k'} \mathrm{Im}(\alpha)] & [10^{k'} \mathrm{Im}(x_n)] & \ldots & [10^{k'} \mathrm{Im}(x_1)] & 0 \\
\end{pmatrix}$$
A similar heuristic argument as before applies. Similarly, if the embedding of the field is real, the last row can be ommitted. 

When the degree of the computed polynomials is small, it is straightforward to verify their correctness; we simply compute the roots of the minimal polynomials and verify that the corresponding points are solutions of the starting system. 

\subsection{Precise Expressions and Continued Fractions} When $3$-torsion points are not defined over number fields of small degree, computing the minimal polynomials of the coefficients is extremely inefficient. However, we note that the 3-torsion points approximated using the method described in the previous subsection 
form a subset which is closed under the action of Galois. Thus, the polynomial whose roots are, say the first coefficients, of the cubics defining the $3$-torsion points in this Galois stable set, has rational coefficients. We approximate these coefficients as complex numbers, using our previously computed high precision approximations. To determine them algebraically, we use the classical theory of continued fractions, which we briefly review below. 

 Let $y_1, ..., y_d$ be complex approximations to the roots of a rational monic polynomial of degree $d$. We can then compute the polynomial $\prod_i (x - y_i)$. This polynomial will be an approximation to the original polynomial, and so it suffices to recognise the coefficients as rational numbers.

 This method is standard and in many undergraduate lecture notes, including \cite{wstein}. Given a real number, $x$, the ``best'' rational approximation comes from the continued fraction expansion. In this case, $x$ is a close approximation to a rational number, $q$. The continued fraction for $q$ terminates, and so the continued fraction expansion for $x$ will, at some point, have a very small error, corresponding to the end point of the continued fraction. We loop through the continued fraction process until an error smaller than a given tolerance is reached, or until a maximum number of iterations is reached. 

 Unfortunately, it is very hard to check this method produces the correct polynomial as they are typically of extremely large degree. There are some heuristics that can provide evidence, for example, checking whether the number of roots of the polynomial over $\mathbb{F}_p$ equals the order of 3-torsion over $\mathbb{F}_p$.
 
\section{Examples} \label{egs}
\subsection{The $3$-torsion Subgroup of the Fermat Quartic}
The Fermat quartic, or equivalently, the modular curve $X_{0}\left( 64 \right)$, has a classical model given by  
\begin{center}
    $x^4 + y^4 - z^4 =0$.
\end{center} 
The transformation $\left( x,y,z \right) \mapsto \left( y, x, y-z \right)$ followed by $\left( x , y , z \right) \mapsto \left( x, \frac{1}{2} y -z, -2z \right)$ moves the rational point $(1:0:1)$ to $P_{0} = (0:1:0)$ with tangent line $z=0$, and we obtain the following model of the curve
\begin{center}
    $ x^{4} - y^{3}z - 4y z^3 = 0. $
\end{center}
We work on the affine patch $z=1$. In this example we were able to use lattice reduction to compute the minimal polynomials of the coefficients of the required cubics.  We found $3$ Galois orbits of $3$-torsion points. First, those cut out by cubics of the form 
\begin{center}
    $ -y^2 + uxy + (72u^6 - 14u^2)y + (1/864)(u^7 - 36u^3)x^3$
\end{center}
where $u^8 -(1/6)u^4 - 1/432 =0 $.  The second orbit consists of the cubics of the form 
\begin{center}
    $-y^2 + uxy + s_{1}(u)y + s_{2}(u)x^3 + s_{3}(u)x^2 + s_{4}(u)x + s_{5}(u)$
\end{center}
where $u$ is a solution of $ x^8 - 8x^7 + 28x^6 - 56x^5 + 52x^4 + 16x^3 - 80x^2 + 64x - 44 =0 $ and $s_{1}, \ldots, s_{5}$ are explicitly given by 
\begin{align*}
& s_{1}(u) = -2u;  \\
 & s_{2}(u) = (1/108)(-u^7 + 7u^6 - 21u^5 + 35u^4 - 26u^3 - 6u^2 - 34u + 46); \ \ \ \ \ \ \\
 & s_{3}(u)  = (1/18)(u^7 - 7u^6 + 21u^5 - 35u^4 + 26u^3 + 6u^2 - 2u + 26); \\ 
 & s_{4}(u) = (1/9)(-u^7 + 7u^6 - 21u^5 + 35u^4 - 26u^3 - 6u^2 + 20u - 26);\\
 &s_{5}(u) = (1/27)(2u^7 - 14u^6 + 42u^5 - 70u^4 + 52u^3 + 12u^2 - 40u + 16).
\end{align*}
The final orbit corresponds to cubics of the form 
\begin{center}
   $ -xyz + u x^{2}z + ( 72u^6 - 14u^2)z^3$
   \end{center}
   and $u$ solves $f = x^8 - (1/6)x^4 -1/432$.
   We verify that the three orbits above generate the entire subgroup $J[3] \cong \left( \mathbb{Z} / 3 \mathbb{Z} \right)^{6}$ by reduction, using the classical injectivity of torsion, see \cite[Appendix]{katz}.

\subsection{The $3$-torsion subgroup of the Klein Quartic} 
We work with the classical affine model of the Klein quartic 
\begin{center}
    $X: x^{3}y + y^{3}z + xz^3 = 0$. 
\end{center}
and compute $\alpha_{1}, \ldots, \alpha_{6}$ such that the cubic defined by 
\begin{center}
    $h = -y^{2} + \alpha_{1}xy + \alpha_{2}y + \alpha_{3}x^{3} + \alpha_{4}x^{2} + \alpha_{5}x + \alpha_{6}$
\end{center}
gives a $3$-torsion element. We computed two orbits of $3-$torsion points, which were sufficient to generate the entire group. The minimal polynomial of $\alpha_{1}$ and  expressions for $\alpha_{2}, \ldots, \alpha_{6}$ in terms of $\alpha_{1}$ can be found in the \texttt{GitHub} repository. 
Both of our orbits are defined over the degree $48$ number field given by 
{\small \begin{center}
     $u^{48} - 32u^{47} + 492u^{46} - 4838u^{45} + 34136u^{44} - 184200u^{43} + 795782u^{42} - 2871488u^{41} + 9027972u^{40} - 25406816u^{39} + 62125865u^{38} - 110851470u^{37} + 20204479u^{36} + 896372004u^{35} - 4721050815u^{34} + 16331796445u^{33} - 43395256901u^{32} + 78805990590u^{31} - 9062375006u^{30} - 597966047368u^{29} + 2489028250704u^{28} - 5314164476524u^{27} + 4247868234838u^{26} + 10286452763613u^{25} - 38984213410246u^{24} + 48194438473527u^{23} + 34337356474920u^{22} - 258320740169787u^{21} + 493721560145034u^{20} - 279315253980219u^{19} - 720711950420342u^{18} + 1333301851341148u^{17} +  650839726478139u^{16} - 3536530808627111u^{15} + 1476853576240844u^{14} + 2773228527944145u^{13} + 2524087738848453u^{12} - 10782439856962015u^{11} + 1577760701045559u^{10} + 10769777804876901u^9 + 3034475537610752u^8 - 18205883402537385u^7 +  351890944113304u^6 + 18574766492209602u^5 - 4386446310081687u^4 - 16132297574659372u^3 + 15821424127424507u^2 - 6053794543417854u + 920504949783049 ;
$
 \end{center}
}
 and as in the previous example, we reduce modulo a prime of good reduction to verify that the two orbits listed above generate the entire $3$-torsion subgroup $J[3] \cong \left( \mathbb{Z} / 3 \mathbb{Z} \right)^{6}$.

\subsection{Generic Examples}
We give two further examples of plane quartic curves with trivial automorphism group where we compute the $3$-torsion polynomial. The curves considered are some of the first in Sutherland's database \cite{sutherland2019database}, with even discriminant and trivial automorphism group: 
\begin{align*}
    & C_{1} :  f_{1} = x^3z+x^2y^2+x^2z^2+xy^3+xyz^2+y^3z+yz^3; \\
    & C_{2} : f_{2} = x^3z+x^2y^2-x^2z^2 + xy^2z+xyz^2+y^4-2y^2z^2-yz^3.    \end{align*}
We computed the $3$-torsion polynomials using the complex approximations and continued fractions method described in the previous section. Notably, in the first example using homotopy continuation we were able to find $726$ approximate complex solutions; and this was sufficient to recover the local Galois representation at $2$, using the method described in the succeeding section. In the second example, we found $728$ solutions and proceeded as expected. See the \texttt{GitHub} repository for our calculations.

% verification 

\section{Local Galois Representations at $3$} \label{galrep}
For a very general plane quartic curve it is impractical to compute explicit generators for the group of three-torsion points, mainly due to the large degree of the scheme. However, the polynomial computed using the method of Section \ref{3torscomps} suffices to recover the local Galois representation at $2$. 

For the rest of this section, we assume that $C$ is a plane quartic over $\Q$ and we assume that $C(\Q) \ne \emptyset$ as in Section \ref{3tors}. The method described in Section \ref{3torscomps} computes a polynomial $F$ of some degree $d$, whose points correspond to $3$-torsion points. We view $C$ as a curve over $\Q_{\ell}$, and note that there is a natural action of $\Gal( \bar{\Q}_{\ell}/ 
\Q_{\ell})$ on the roots of $F$, which corresponds to the Galois action on $J[3]$. In this section, we explain how one can recover the Galois representation $\rho: \Gal(\bar{\Q}_{\ell} / \Q_{\ell}) \rightarrow \Aut(J[3])$ for any $\ell \ne 2$, from the Galois action on the roots of $F$. We use this local representation to determine the wild conductor exponent away from $3$. 
\begin{remark}
Generically, $d = 728$, however this is not always the case as previously explained, and in fact, for our calculations $d \ge 243$ suffices.
\end{remark}

 \subsection{Splitting Field and Galois Group}
Given $F$, we compute its splitting field $K$ over $\Q_{\ell}$, the set of roots $R$ over $K$ and the Galois group $G$ of $F$ acting on $R$. In \texttt{Magma}, this is done using the algorithms given in \cite{ford2002fast} and \cite{pauli2001factoring}. In some of our examples, we found it more efficient to compute the splitting field directly. Our implementation of this is given in the repository, where we compute these factorisations via Hensel lifting and Newton polygons. Finding all the roots of $F$ using Hensel lifting can be intensive, as many start solutions are needed over highly ramified extensions. We exploit the factorisable structure to compute one for each irreducible factor, and then compute its Galois orbit.

The Galois group is computed from the splitting field. Finally, we must compute the permutation action of the Galois group on the roots. The naive approach of computing this requires computing the differences between the image of each root under Galois, and all the other roots. Arithmetic in large field extensions is slow, so we use their images in the residue field to quickly perform comparisons.

\subsection{Galois groups and $\text{GSp}_6(\mathbb{F}_{3})$}\label{grp_ident}
 If the degree of $F$ is at least $3^5$ and $F$ is square free, then the three-torsion points corresponding to the roots of $F$ generate the entire $3$-torsion subgroup and computing the Galois action on the roots is sufficient to determine the Galois representation. This is since there are more roots of $F$ than non-zero vectors in the largest proper subspace of $\mathbb{F}_3^6$, and so the torsion points corresponding to these roots must span the whole space. 

The Galois representation attached to $J[3]$ naturally has image in $\text{GSp}_6(\mathbb{F}_3)$, as it preserves the Weil pairing. We therefore wish to identify the Galois group, which is computed as a permutation group, as a subgroup of this. As our computations record little-to-none of the vector space structure, it certainly does not contain a choice of basis, and so we can only identify it with subgroups of $\text{GSp}_6(\mathbb{F}_3)$ up to conjugation by $\text{GL}_6(\mathbb{F}_3)$. Further, since the absolute Galois group of $\mathbb{Q}_\ell$ is soluble, we can restrict to the soluble subgroups of $\text{GSp}_6(\mathbb{F}_3)$. Our approach is to pre-compute all these subgroup classes, along with certain invariants, so that we can compare against the computed invariants of the Galois group. We save these subgroups separated by order, as there are so many (approximately 16000) that loading all of them at once is prohibitively slow.

For each equivalence class of subgroups under conjugation, $H \leq \text{GSp}_6(\mathbb{F}_3)$, we compute the order of $H$, the sizes of its orbits on $\mathbb{F}_3^6 \setminus \{0\}$, and the permutation representation on this set. Given a permutation group on $n \leq 728$ symbols, we first restrict to subgroups of the same order, which as remarked above, are stored separately. From this table of subgroups, we compare the orbit sizes: for every orbit of the permutation group, there must exist an orbit of $H$ of the same size if $H$ is the image of the local Galois representation. Finally, for those subgroups satisfying this, we further restrict to those abstractly isomorphic to the permutation group. These isomorphisms are not necessarily unique, for example, by composing with conjugation by a fixed element. In fact, any two isomorphisms disagree by an automorphism of the group. We compute representatives for the outer automorphism group of the permutation group (the group of automorphisms modulo the subgroup of conjugations), as these will give all the isomorphisms where the action of individual elements may differ. For each isomorphism, we compare the orbit sizes for each conjugacy class of permutations, and its image. As before, if this combination of subgroup and isomorphism realise the permutation group as a matrix group, then for every cycle of the permutation, there exists a cycle of the same size in the permutation action of its image in $H$. 

In practice, this identifies the subgroup uniquely, and there are only a few choices for the isomorphism, which do not effect the computed conductor (see later). This ambiguity is likely due to conjugation by elements of $\text{GSp}_6(\mathbb{F}_3)$ or $\text{GL}_6(\mathbb{F}_3)$, which are not accounted for when computing the automorphism groups.

\subsection{Inertia and Ramification Groups}
We recall the definitions of the inertia and ramification groups, and how to calculate them in \texttt{Magma}.

The inertia subgroup, $G_0$, is the subgroup of $G = \mathrm{Gal}(\mathbb{Q}_{\ell}(J[3]) / \mathbb{Q}_{\ell})$ consisting of elements, $\sigma$, such that $v_\ell(x - \sigma(x)) \geq 1$ for all $x$ in the ring of integers  of $ \mathbb{Q}_{\ell}(J[3])$. Its fixed field is the maximal unramified subfield of $\mathbb{Q}_\ell(J[3])$. Similarly, we define the ramification groups in the lower numbering, 
\begin{center}
    $G_i = \{ \sigma \in G \ \vert v_\ell(x - \sigma(x)) \geq i + 1 \ \text{for all} \  x \in \mathcal{O}_{\mathbb{Q}_{\ell}(J[3])} \}  $
\end{center}
for all  for all non-negative integers $i$. 
For any real number $u \ge -1$, set $G_u = G_{\lceil u \rceil}$ and define $\phi : [-1, \infty) \rightarrow [-1, \infty )$ by
\begin{center}
    $\phi (u) = \displaystyle \int_{0}^{u} \frac{1 }{[G_0 : G_t]} dt$ 
\end{center}
if $u \ge 0$ and $\phi (u) =u$ if $-1 \le u \le 0$. The ramification groups in upper numbering are defined as 
\begin{center}
    $G^i = G_{\phi(i)}$
\end{center}
for all non-negative integers $i$.

The computer algebra package $\texttt{Magma}$ can compute ramification groups, however due to the work-arounds mentioned in the section on computing splitting fields, this option is not available to us. Instead, we compute these subgroups using a variation of the above definitions. For an element, $\sigma$, of the Galois group to belong to $G_i$ it is enough to satisfy $v_\ell(\alpha - \sigma(\alpha)) \geq i + 1$, where $\alpha$ generates ring of integers of the splitting field. We use this combined with the isomorphism identified in the previous section to determine the ramification groups and their action on the 3-torsion.

\section{Two-Torsion Subgroup and the Local Galois Representation at $2$}
 \label{2torsion}
 For plane quartics, $2$-torsion points have a beautiful geometric description, and they are significantly more manageable to compute. We'll review this here, and with a view towards computing $n_{p, \text{wild}}$ for all odd primes $p$ in the sections which follows.  A detailed overview of bitangents and their connection to $2$-torsion can be found in  \cite[Chapter 5]{dog} and \cite{caporaso}, and computational examples can be found in \cite{lupoiantwo}.
 \subsection{Bitangents and $2$-torsion}
Let $l$ and $k$ be bitangents to the curve, that is lines in $\mathbb{P}^{2}$ intersecting the curve at $2$ points, each with multiplicity $2$. The function on the curve defined by the quotient of the two lines is of the form 
\begin{center}
    $ \divv \left( \frac{l}{k} \right) = 2 (P_{1} + P_{2} - Q_{1} - Q_{2}) = 2D_{l,k} $
\end{center}
for some (not necessarily distinct) points $P_{i},Q_{j} \in C$. It is clear that $D_{l,k}:= P_1 + P_2 -Q_1 - Q_2$ defines a $2-$torsion point on the Jacobian. In fact, such points generate the entire $2$-torsion subgroup. 
\begin{theorem}
 The two torsion subgroup of the Jacobian of a plane quartic is generated by divisors of the form $D_{l,k}$, where $l$ and $k$ are bitangents to the curve.  
\end{theorem}
   \begin{proof}
Bitangent lines correspond to odd theta characteristics, that is, square roots of the canonical divisor, whose Riemann-Roch space has odd dimension. As in the construction above, one can see that the difference of two such divisors defines a two torsion point on the Jacobian. The fact that such differences generate the entire two-torsion subgroup is a classical result, proved for instance in \cite[Chapter 5]{dog}.
\end{proof} 

To compute the bitangents to $C$ we proceed as before. We first derive a system of equations whose solutions correspond to bitangents to the curve. Suppose that $a_{1}x + a_{2}y + a_{3}z$ defines a bitangent to the curve. We work on the affine patch $z=1$, and suppose that $a_{1} \ne 0$. Then rescaling if necessary, we search for $\alpha_{1}, \alpha_{2} \in \Bar{\Q}$ such that 
\begin{center}
    $f(\alpha_{1}y + \alpha_{2},y,1) = \alpha_{3}(x^2 + \alpha_{4}x + \alpha_{5})^{2} $
\end{center}
for some $\alpha_{3}, \alpha_{4} , \alpha_{5} \in \Bar{\Q}$.
Equating coefficients in the above expression, we obtain a system of five equations in the $\alpha_{1}, \ldots \alpha_{5}$, with $-x + \alpha_{1}y + \alpha_{2}$ defining a bitangent to the curve. We compute a Groebner basis of the ideal defined by these equations (using \texttt{Magma}'s implementation of \texttt{GroebnerBasis()}). Note that if one does not obtain all $28$ bitangents via the above system, the above can be repeated working on a different affine patch, or by assuming a different coefficient to be non-zero, or by composing $f$ with a generic element of $\text{PSL}_{2}(\mathbb{Z})$, to obtain an isomorphic curve. In practice, it is this last option that is most viable as a random element of $\mathrm{PSL}_2(\mathbb{Z})$ can be chosen repeatedly until a good choice of model is found.

\subsection{Local Galois Representation at $2$}
To compute $J[2]( \Bar{\mathbb{Q}})$, one can work with this Groebner basis, compute the relevant splitting field and determine the bitangents explicitly. Although this is computationally feasible, we found it to be very inefficient for our computations, and thus to determine $n_{\text{wild},p}$ we aimed to localise as many of our computations as possible. 

From the Groebner basis (computed over $\mathbb{Q}$), our bitangents are of the form 
\begin{center}
    $-x + f_{1}(u)y + f_{2}(u)z $
\end{center}
where $u$ is the solutions of a degree $28$ polynomial $f_{3}(x)$, with $f_{i}(x) \in \mathbb{Q}[x]$. We compute the splitting of $f_{3}$ over $\Q_{p}$ and the corresponding bitangents. To determine the action of $\Gal( \Bar{\mathbb{Q}}_{p} / \mathbb{Q}_{p})$ on $J[2]$ efficiently, we first compute a basis for $J[2]$ using the bitangents. We fix a bitangent, $l_0$, and from the sequence $D_{l_1, l_0}, ..., D_{l_m, l_0}$, where the $l_i$ run over bitangents not equal to $l_0$, we construct a basis as follows. We add the first two divisors to the basis set (as they are non-zero and distinct, they must be independent over $\mathbb{F}_2$). Then, for each subsequent element of this list, we test if it is the sum of any subset from the basis so far. If it can be represented as such as sum, we store this representation, and otherwise, we add the divisor to the basis. The key to computing this linear equivalence is the following lemma.

\begin{proposition}
    The sum $D_{k_1, l_0} + ... + D_{k_n, l_0}$ is equal to $D_{l, l_0}$ if and only if there exists a curve (not containing $C$ as an irreducible component) through the intersection points of $C$ with $k_1$, ..., $k_n$, and $l$ of degree $\frac{n + 1}{2}$ if $n$ is odd, and through $k_1$, ..., $k_n$, $l$ and $l_0$ of degree $\frac{n + 2}{2}$ if $n$ is even.
\end{proposition}
\begin{proof}
    If $n$ is odd, then the equality $D_{k_1, l_0} + ... + D_{k_n, l_0} \sim D_{l, l_0}$ can be rearranged to $D_{k_1, l_0} + ... + D_{k_n, l_0} + D_{l, l_0} \sim 0$. Let $D'_{l}$ denote the degree 2 divisor giving the intersection points of a bitangent $l$ with $C$. The desired linear equivalence is equivalent to $D'_{k_1} + ... + D'_{k_n} + D'_{l} \sim (n + 1) D'_{l_0}$. As $l_0$ is a bitangent, the right hand side is equivalent to $\frac{n + 1}{2} K$. Finally, the linear system of $\frac{n + 1}{2} K$ is precisely the intersection of the quartic with degree $\frac{n + 1}{2}$ curves (not containing $C$), the result follows. 

    The argument for even $n$ is essentially identical for odd $n$, except the equality $D'_{k_1} + ... + D'_{k_n} + D'_{l} \sim (n + 1)D'_{l_0}$ is replaced with $D'_{k_1} + ... + D'_{k_n} + D'_{l} + D'_{l_0} \sim (n + 2)D'_{l_0} \sim \frac{(n + 2}{2}K$.
\end{proof}

This reduces linear equivalence to a problem of curve fitting. Curve fitting can be reduced to linear algebra over the splitting field as follows. Since the basis has size at most 6, we will be fitting at worst degree 4 curves. For each bitangent, we compute the intersection points with $C$, and the value of each monomial of degree at most 4 is computed for each. For a hyperflex, where there is a unique intersection point, we also record the derivatives of these monomials along the tangent direction. Thus, to each bitangent we obtain two linear equations which the coefficients of a curve of degree at most 4 passing through its intersection points must satisfy. When the intersection points are not defined over the splitting field it is desirable to combine the two equations to get two equations defined over the smaller field. Determining linear equivalence is now a case of determining the nullity of the corresponding matrix.

The Galois group $\Gal(\mathbb{Q}_{p}(J[2]) / \mathbb{Q}_{p})$ and its ramification groups can be computed using \texttt{Magma's} local fields implementation. The action of these groups on the two torsion can be extracted using the basis and computed expressions found earlier. In particular, the image of each basis element can be computed, and re-expressed in terms of the basis as $\sigma(D_{k, l_0}) = \frac{1}{2} \mathrm{div}(\sigma(k) / \sigma(l_0)) = D_{k', l_0} - D_{l'_0, l_0}$, where $k'$ and $l'_0$ are the indices of $\sigma(k)$ and $\sigma(l_0)$, and $\sigma$ is an element of the Galois group. This gives an embedding of the abstract Galois group into $\mathrm{GL}_6(\mathbb{F}_2)$, and the ramification groups can be pushed forward under this inclusion.

\begin{remark}
In principle, combining all the local Galois actions (along with calculations of the Weil pairing) with the Galois group of the degree 28 polynomial and applying a method similar to that in Section~\ref{grp_ident}, would allow for the identification of the image of the mod-2 Galois representation over $\mathbb{Q}$, without computing the global splitting field.
\end{remark}

\section{Local Conductor Exponents} \label{cond}
The Galois action on the $2$ and $3$ torsion subgroup of $J$ are sufficient to completely determine the wild conductor exponent at any prime as we now explain. 

Let $C$ be any smooth, projective curve of genus $g$ defined over $\mathbb{Q}$.  The conductor of $C$ is an arithmetic invariant of the curve, defined as a product over the primes of bad reduction for the curve. The problem of computing conductors is essentially the problem of computing the local conductor exponent $n_{p}$ at each prime of bad reduction for the curve. In this section, we give an overview of relevant definitions and results, see \cite{ulmer} for details. 

To compute $n_{p}$, we view $C$ as a smooth projective curve over $\mathbb{Q}_{p}$, let $J$ be the Jacobian variety associated to $C$ and for a prime $\ell \ne p$ let $T = T_{\ell}J$ and $V = T \otimes_{\mathbb{Z}_{\ell}} \mathbb{Q}_{\ell}$ be the associated $\ell-$adic Tate module and $\ell-$adic representation respectively. Then $n_{p}$ is the Artin conductor of $V$ defined as follows
\begin{center}
$n_p = \displaystyle \int_{-1}^{\infty} \text{codim}  V^{G_{\mathbb{Q}_p}^{u}} \ du$
\end{center}
where $G_{\mathbb{Q}_p} = \Gal \left( \overline{\mathbb{Q}_p} / \mathbb{Q}_p \right)$ is the absolute Galois group of $K$ and $ \lbrace G_{\mathbb{Q}_p}^{u} \rbrace_{u \ge -1}$ denote the ramification groups of $G_{\mathbb{Q}_p}$ in upper numbering.
The definition is independent of the choice of prime $\ell$.

\subsection{Tame Conductor Exponent}
For computational purposes, it's convenient to break up the above quantity into two parts. First, we define the `tame part' of $n_{p}$ as 
\begin{center}
$n_{p, \text{tame}} =  \displaystyle \int_{-1}^{0} \text{codim}  V^{G_{\mathbb{Q}_p}^{u}} \ du $.
\end{center}
We recall that $G_{\mathbb{Q}_p}^{0} = I $ is the inertia subgroup, and hence the tame conductor exponent is simply 
\begin{center}
  $n_{p, \text{tame}} = 2g - \text{dim}V^{I}$.
\end{center}
This is computable directly from a simple normal crossings regular model of the curve over $\mathbb{Z}_{\ell}$, as we can deduce the following invariants:
\begin{itemize}
\item the abelian part $a$,  equal to the sum of the genera of all components of the model,
\item the toric part $t$,  equal to the number of loops in the dual graph of the regular model of $C$.
\end{itemize}
The tame part of the exponent is equal to $2g -2a -t$, see \cite[Chapter 9]{NeronModels} (for the relationship between the abelian and toric ranks and the special fibre) and \cite{Lang} (for the tame conductor formula, which there is given in terms of the unipotent dimension, $g - a - t$). Regular models can be computed in principle by taking any model of the curve and performing repeated blowups until it becomes regular. A number of works have made significant progress towards making regular models accessible, see for instance \cite{ModelsoverDVRs}. However, to date, implementations are slow and inefficient. It is for this reason that we shall not discuss the tame part further. 
\subsection{Wild Conductor Exponent}
The second part of the conductor exponent is the `wild part', defined as
\begin{center}
    $n_{\text{wild}}  =  \displaystyle \int_{0}^{\infty} \text{codim}  V^{G_{\mathbb{Q}_p}^{u}} \ du$.
\end{center}

As the higher ramification groups are pro-$p$, the dimension of the fixed part of the representation is the same as the dimension of the fixed part of the $\ell$-torsion. We may therefore replace the absolute Galois group with the Galois group of the finite extension $\mathbb{Q}_p(J[\ell])$, and can swap to lower numbering. The full details of this calculation can be found in \cite{ulmer}.
\begin{center}
   $n_{p, \text{wild}} = \displaystyle \sum_{i=1}^{\infty} \frac{\text{codim} J\left[ \ell \right]^{G_{i}}}{[ G_{0} : G_{i}]} $.
\end{center}
Thus for a given curve explicitly knowing the $\mathbb{Q}_{p}$ -Galois action on  $J[ 2]$ and $J[ 3]$ is sufficient to determine all $n_{p, \text{wild}}$. In particular, in the case of plane quartics, one can deduce $n_{l,\text{wild}}$ for all $l \ne 3$, from the $3$-torsion polynomial, as explained in the previous section. In practice we only use the $3$-torsion polynomial to determine $n_{2, \text{wild}}$, and compute all other wild conductor exponent from $J[2]$, as this is particularly efficient to compute in the case of plane quartics. 

\begin{remark}
Since we only use the 3-torsion to compute the wild conductor at 2, we can extend the ground field to any field unramified at 2. This can facilitate computing the conductor for curves with no rational points, as they will often have points over a small field unramified at 2. This will, unfortunately, make recognising the roots as algebraic integers harder. The most straightforward solution is to compute the minimal polynomial as in the continued fraction method, but recognise its coefficients using the lattice reduction techniques.
\end{remark}

\subsection{Examples}
\subsubsection{The Fermat and Klein Quartic}
When the explicit generators of the $3$-torsion subgroup are known, computing the wild conductor exponent at $2$ is a simple calculation. For the Fermat quartic, using \texttt{Magma} we compute the extension $\mathbb{Q}_{2}\subseteq K$ over which the $3$-torsion generators are defined, and further compute $G = \Gal \left( K / \Q_{2} \right)$ and the ramification subgroups  $\{ G_{i} \}_{i \ge -1}$. We found that $G_{n} =1 $ for all $n \ge 8$ and $\vert G_{0}\vert = \vert G_{1} \vert = 16$, $\vert G_{2} \vert = \vert G_{3} \vert =4 $ and $\vert G_{4} \vert = \vert G_{5} \vert = \vert G_{6} \vert = \vert G_{7} \vert =2$. Explicitly taking Galois invariants, we find:
\begin{align*}
& J[3]^{G_{0}} =  J[3]^{G_{2}}  =  \{ 0 \}, \\
& J[3]^{G_{4}} \cong \left( \mathbb{Z} / 3 \mathbb{Z} \right)^{4},
\end{align*}
and thus $ n_{2, \text{wild}} =  6 + 2 \cdot \frac{6}{4} + 4 \cdot  \frac{2}{8} = 10 $. Observe that as $J[3]^{G_0} = \{0\}$, since  the 3-adic Tate module maps onto the 3 torsion, $V^{G_0} = \{0\}$ and so the tame part of the conductor is 6.  Combining this gives the conductor exponent at 2 as $n_2 = 16$. Since the discriminant of the standard model of the Fermat quartic  is $2^{40}$, $2$ is the only prime of bad reduction. 

As for the Klein quartic, the discrimiant of the standard model is $7^7$. We note that number field over which the $3$-torsion is defined is unramified above 2, the ramification groups beyond $G_0$ are trivial. In particular, the wild conductor at 2 is trivial, thus verifying the fact that the curve  has good reduction at 2.  Instead, we demonstrate the method at 7, the only bad prime for this curve. Let $G = \mathrm{Gal}(\mathbb{Q}_7(f) / \mathbb{Q}_7)$, where $f$ is the defining polynomial of $K$ from Section 6. We denote the ramification groups in the lower numbering by $\{G_i\}_{i \geq -1}$. As $G$ has order 24, which is co-prime to 7, $G_1$ is trivial since it is a $7$-group. This shows the wild conductor is again, trivial.

As in the previous example, we compute that $J[3]^{G_0} = \{0\}$, and so there are no non-trivial fixed points of the Tate module. This shows $n_{7, \mathrm{tame}} = 6$. 

\subsubsection{Generic Examples}
For the first 100 entries of Sutherland's database \cite{sutherland2019database}  we computed the wild conductor exponents at all odd primes dividing the discriminant using our implementation. We found that all wild exponents were $0$, except for the few exceptions recorded in the table below. The computations took approximately 4.5 hours. We believe that our implementation can be improved greatly if one could reduce the reliance on Magma's local fields commands. In particular, we found the computation of the splitting field of the bitangents to be very computationally expensive as well as determining a basis, due to a reliance on inexact linear algebra.  We note that the initial model of the plane quartic does affect the performance of our algorithm, mostly due to the fact that the precision required to compute the bitangents depends on the coefficients of the equations. As our implementation involves a random choice of model, there are cases for which it is particularly slow. 

\begin{table}[ht]
    \centering
    \begin{tabular}{|c|c|c| }
    \hline
    Discriminant & $f$ & Exceptional Exponent \\ 
    \hline 
$15957$ & $x^3z + 5x^2yz + 2xy^2z + 8xyz^2 - y^4 - 4y^3z - 8y^2z^2 + yz^3$ & $n_{3} = 2$ \\
    \hline
$15957$ & $ x^3z + x^2yz + x^2z^2 + xy^2z - xyz^2 + xz^3 - y^4 + y^3z$ & $n_{3} = 2$ \\
\hline 
$17307$ &  $x^3z + x^2y^2 + x^2yz + x^2z^2 + xy^3 + xy^2z + 
    xyz^2 + y^3z + y^2z^2 + yz^3$ & $ n_{3} =1$ \\
    \hline 
$ 20331$ & $x^3z + x^2y^2 - x^2z^2 - xy^2z + xz^3 + y^3z + y^2z^2$ & $n_{3} =2$ \\
\hline 
$22707$ & $x^3z + x^2y^2 + x^2yz + x^2z^2 - xy^3 - xy^2z - 
    xyz^2 - y^3z - y^2z^2 - yz^3$ &  $ n_3 = 1$ \\
    \hline 
$25029$ & $x^3z + x^2y^2 + x^2z^2 + 2xyz^2 + y^3z - y^2z^2 + yz^3$ & $n_{3} =3 $ \\
\hline 
$27702$ &  $x^3z + x^2y^2 + x^2yz + x^2z^2 + xy^2z + xyz^2 + xz^3 + y^3z + y^2z^2 + yz^3$ &  $n_{3} =2$ \\
\hline
    \end{tabular}
    \end{table}
    
For the two generic examples discussed in Section \ref{egs}, we compute the $n_{2, \text{wild}}$ and found it to be $2$ in both cases. 

\section*{Declarations}
\subsection*{Conflict of Interest} The authors state that there are no conflicts of interest related to this work.
\subsection*{Data availability} Data sharing is not applicable to this article, there are not datasets generated. The examples discussed in this article and the accompanying code can be found in the online repository 
\begin{center}
\href{https://github.com/ElviraLupoian/Genus3Conductors}{https://github.com/ElviraLupoian/Genus3Conductors}
\end{center}
\bibliographystyle{abbrv}
\bibliography{ref}

\begin{thebibliography}{10}

\bibitem{bernard2009jacobians}
N.~Bernard, F.~Lepr{\'e}vost, and M.~Pohst.
\newblock {Jacobians of genus-2 curves with a rational point of order 11}.
\newblock {\em Experimental Mathematics}, 18(1):65--70, 2009.

\bibitem{NeronModels}
S.~Bosch, W.~L{\"u}tkebohmert, and M.~Raynaud.
\newblock {\em {N{\'e}ron models}}, volume~21.
\newblock Springer Science \& Business Media, 2012.

\bibitem{JulHC}
P.~Breiding and S.~Timme.
\newblock {{H}omotopy{C}ontinuation.jl: {A} {P}ackage for {H}omotopy
  {C}ontinuation in {J}ulia}.
\newblock In {\em International Congress on Mathematical Software}, pages
  458--465. Springer, 2018.

\bibitem{bruin2014descent}
N.~Bruin, E.~V. Flynn, and D.~Testa.
\newblock {Descent via (3, 3)-isogeny on Jacobians of genus 2 curves}.
\newblock {\em Acta Arithmetica}, 165:201--223, 2014.

\bibitem{BS2tors}
N.~Bruin and M.~Stoll.
\newblock {Two-cover descent on hyperelliptic curves}.
\newblock {\em Math. Comp.}, 78(268):2347--2370, 2009.

\bibitem{caporaso}
L.~Caporaso.
\newblock {On modular properties of odd theta-characteristics}.
\newblock In {\em Advances in algebraic geometry motivated by physics
  ({L}owell, {MA}, 2000)}, volume 276 of {\em Contemp. Math.}, pages 101--114.
  Amer. Math. Soc., Providence, RI, 2001.

\bibitem{ModelsoverDVRs}
T.~Dokchitser.
\newblock Models of curves over discrete valuation rings.
\newblock {\em Duke Mathematical Journal}, 170(11):2519--2574, 2021.

\bibitem{ddmm}
T.~Dokchitser, V.~Dokchitser, C.~Maistret, and A.~Morgan.
\newblock Arithmetic of hyperelliptic curves over local fields.
\newblock {\em Mathematische Annalen}, 385(3):1213--1322, 2023.

\bibitem{dokdor}
T.~Dokchitser and C.~Doris.
\newblock {3-torsion and conductor of genus 2 curves}.
\newblock {\em Mathematics of Computation}, 88(318):1913--1927, 2019.

\bibitem{dog}
I.~V. Dolgachev.
\newblock {\em {Classical algebraic geometry}}.
\newblock Cambridge University Press, Cambridge, 2012.
\newblock A modern view.

\bibitem{flynn1990large}
E.~V. Flynn.
\newblock {Large rational torsion on abelian varieties}.
\newblock {\em Journal of Number Theory}, 36(3):257--265, 1990.

\bibitem{flynn2015descent}
E.~V. Flynn.
\newblock {Descent via (5, 5)-isogeny on Jacobians of genus 2 curves}.
\newblock {\em Journal of Number Theory}, 153:270--282, 2015.

\bibitem{ford2002fast}
D.~Ford, S.~Pauli, and X.-F. Roblot.
\newblock A fast algorithm for polynomial factorization over $\mathbb{Q}_p$.
\newblock {\em Journal de th{\'e}orie des nombres de Bordeaux}, 14(1):151--169,
  2002.

\bibitem{howe2015genus}
E.~W. Howe.
\newblock {Genus-2 Jacobians with torsion points of large order}.
\newblock {\em Bulletin of the London Mathematical Society}, 47(1):127--135,
  2015.

\bibitem{howe2000large}
E.~W. Howe, F.~Lepr{\'e}vost, and B.~Poonen.
\newblock {Large torsion subgroups of split Jacobians of curves of genus two or
  three}.
\newblock {\em Forum Mathematicum}, 12(3):315--364, 2000.

\bibitem{katz}
N.~M. Katz.
\newblock {Galois properties of torsion points on abelian varieties}.
\newblock {\em Invent. Math.}, 62(3):481--502, 1981.

\bibitem{Lang}
S.~Lang.
\newblock Abelian varieties.
\newblock In {\em Number Theory III: Diophantine Geometry}, pages 68--100.
  Springer Berlin Heidelberg, Berlin, Heidelberg, 1991.

\bibitem{leprevost1992torsion}
F.~Lepr{\'e}vost.
\newblock {Torsion sur des familles de courbes de genre g}.
\newblock {\em Manuscripta mathematica}, 75:303--326, 1992.

\bibitem{lupoiantwo}
E.~Lupoian.
\newblock {Two-torsion subgroups of some modular Jacobians}.
\newblock {\em International Journal of Number Theory}, 20(10):2543--2573,
  2024.

\bibitem{lupoianthree}
E.~Lupoian.
\newblock Three-torsion subgroups and wild conductors of genus 3 hyperelliptic
  curves.
\newblock {\em Journal of Number Theory}, 278:267--284, 2026.

\bibitem{mazur1977modular}
B.~Mazur.
\newblock {Modular curves and the Eisenstein ideal}.
\newblock {\em Publications Math{\'e}matiques de l'Institut des Hautes
  {\'E}tudes Scientifiques}, 47(1):33--186, 1977.

\bibitem{j2023computing}
J.~S. M{\"u}ller and B.~Reitsma.
\newblock Computing torsion subgroups of jacobians of hyperelliptic curves of
  genus 3.
\newblock {\em Research in Number Theory}, 9(2):23, 2023.

\bibitem{pauli2001factoring}
S.~Pauli.
\newblock Factoring polynomials over local fields.
\newblock {\em Journal of Symbolic Computation}, 32(5):533--547, 2001.

\bibitem{silv2}
J.~H. Silverman.
\newblock {\em {Advanced topics in the arithmetic of elliptic curves}}, volume
  151.
\newblock Springer Science \& Business Media, 1994.

\bibitem{smart}
N.~P. Smart.
\newblock {\em {The algorithmic resolution of {D}iophantine equations}},
  volume~41 of {\em London Mathematical Society Student Texts}.
\newblock Cambridge University Press, Cambridge, 1998.

\bibitem{wstein}
W.~A. Stein.
\newblock {Recognizing Rational Numbers From Their Decimal Expansion}.
\newblock \url{https://wstein.org/edu/2007/spring/ent/ent-html/node74.html},
  2007.
\newblock Accessed 22/01/24.

\bibitem{stoerNA}
J.~Stoer, R.~Bulirsch, R.~Bartels, W.~Gautschi, and C.~Witzgall.
\newblock {\em {Introduction to numerical analysis}}, volume 1993.
\newblock Springer, 1980.

\bibitem{stoll1998height}
M.~Stoll.
\newblock {On the height constant for curves of genus two}.
\newblock {\em Acta Arithmetica}, 90:183--201, 1998.

\bibitem{sutherland2019database}
A.~V. Sutherland.
\newblock {A database of nonhyperelliptic genus 3 curves over Q}.
\newblock In {\em Proceedings of the Thirteenth Algorithmic Number Theory
  Symposium}, volume~2, pages 443--459. Mathematical Sciences Publishers
  Berkeley, CA, 2019.

\bibitem{ulmer}
D.~Ulmer.
\newblock {Conductors of l-adic representations}.
\newblock {\em Proceedings of the American Mathematical Society},
  144:2291--2299, 2016.

\bibitem{HC}
J.~Verchelde.
\newblock Homotopy continuation methods for solving polynomial systems.
\newblock {\em PhD thesis}, 1996.

\end{thebibliography}
\end{document}